\documentstyle[12pt, amstex, amssymb, amscd]{amsart}

\newtheorem{theorem}{Theorem}[section]
\newtheorem{Theorem}[theorem]{Theorem}
\newtheorem{Problem}[theorem]{Problem}
\newcommand{\BProblem}{\begin{Problem}}
\newcommand{\EProblem}{\end{Problem}}
\newtheorem{Vermutung}[theorem]{Vermutung}
\newtheorem{Folgerung}[theorem]{Folgerung}
\newtheorem{Conjecture}[theorem]{Conjecture}
\newtheorem{Proposition}[theorem]{Proposition}
\newtheorem{Corollary}[theorem]{Corollary}
\newtheorem{Korollar}[theorem]{Korollar}
\newtheorem{Lemma}[theorem]{Lemma}
\newtheorem{Satz}[theorem]{Satz}
\newcommand{\BSatz}{\begin{Satz}}
\newcommand{\ESatz}{\end{Satz}}

\newtheorem{Hauptsatz}[theorem]{Hauptsatz}
\newcommand{\BHauptsatz}{\begin{Hauptsatz}}
\newcommand{\EHauptsatz}{\end{Hauptsatz}}

\newtheorem{Klassifikationssatz}[theorem]{Klassifikationssatz}
\newcommand{\BKlassifikationssatz}{\begin{Klassifikationssatz}}
\newcommand{\EKlassifikationssatz}{\end{Klassifikationssatz}}
\newtheorem{Hilfssatz}[theorem]{Hilfssatz}
\newcommand{\BHilfssatz}{\begin{Hilfssatz}}
\newcommand{\EHilfssatz}{\end{Hilfssatz}}

\newtheorem{Corollary (of the proof)}[theorem]{Corollary (of the proof)}

\theoremstyle{definition}
\newtheorem{Definition}[theorem]{Definition}

\theoremstyle{remark}
\newtheorem{Remark}[theorem]{Remark}

\newcommand{\ERemark}{\end{Remark}}
\newcommand{\BRemark}{\begin{Remark}}
\newtheorem{Remarks}[theorem]{Remarks}
\newcommand{\ERemarks}{\end{Remarks}}
\newcommand{\BRemarks}{\begin{Remarks}}

\newtheorem{Ramifications}[theorem]{Ramifications}
\newcommand{\BRamifications}{\begin{Ramifications}}
\newcommand{\ERamifications}{\end{Ramifications}}
\newtheorem{Bemerkung}[theorem]{Bemerkung}
\newcommand{\BBemerkung}{\begin{Bemerkung}}
\newcommand{\EBemerkung}{\end{Bemerkung}}

\newtheorem{Anschauung}[theorem]{Anschauung}
\newcommand{\BAnschauung}{\begin{Anschauung}}
\newcommand{\EAnschauung}{\end{Anschauung}}

\newtheorem{Bemerkungen}[theorem]{Bemerkungen}
\newcommand{\BBemerkungen}{\begin{Bemerkungen}}
\newcommand{\EBemerkungen}{\end{Bemerkungen}}
\newtheorem{Kommentar}[theorem]{Kommentar}
\newcommand{\BKommentar}{\begin{Kommentar}}
\newcommand{\EKommentar}{\end{Kommentar}}
\newtheorem{Notation}[theorem]{Notation}
\newcommand{\BNotation}{\begin{Notation}}
\newcommand{\ENotation}{\end{Notation}}
\newtheorem{Behauptung}[theorem]{Behauptung}
\newcommand{\BBehauptung}{\begin{Behauptung}}
\newcommand{\EBehauptung}{\end{Behauptung}}
\newtheorem{Beispiel}[theorem]{Beispiel}
\newcommand{\BBeispiel}{\begin{Beispiel}}
\newcommand{\EBeispiel}{\end{Beispiel}}

\newtheorem{Motivation}[theorem]{Motivation}
\newcommand{\BMotivation}{\begin{Motivation}}
\newcommand{\EMotivation}{\end{Motivation}}

\newtheorem{Beispiele}[theorem]{Beispiele}
\newcommand{\BBeispiele}{\begin{Beispiele}}
\newcommand{\EBeispiele}{\end{Beispiele}}
\newtheorem{Example}[theorem]{Example}
\newcommand{\EExample}{\end{Example}}
\newcommand{\BExample}{\begin{Example}}
\newtheorem{Examples}[theorem]{Examples}
\newcommand{\EExamples}{\end{Examples}}
\newcommand{\BExamples}{\begin{Examples}}
\newtheorem{Exercise}[theorem]{Exercise}
\newcommand{\EExercise}{\end{Exercise}}
\newcommand{\BExercise}{\begin{Exercise}}
\newtheorem{Ubung}[theorem]{"Ubung}
\newcommand{\EUbung}{\end{Ubung}}
\newcommand{\BUbung}{\begin{Ubung}}
\newtheorem{Ubungen}[theorem]{"Ubungen}
\newcommand{\EUbungen}{\end{Ubungen}}
\newcommand{\BUbungen}{\begin{Ubungen}}
\newtheorem{Exercises}[theorem]{Exercises}
\newcommand{\EExercises}{\end{Exercises}}
\newcommand{\BExercises}{\begin{Exercises}}
\newtheorem{Claim}[theorem]{Claim}

\newcommand{\BTheorem}{\begin{Theorem}}
\newcommand{\ETheorem}{\end{Theorem}}
\newcommand{\BVermutung}{\begin{Vermutung}}
\newcommand{\EVermutung}{\end{Vermutung}}
\newcommand{\BFolgerung}{\begin{Folgerung}}
\newcommand{\EFolgerung}{\end{Folgerung}}
\newcommand{\BConjecture}{\begin{Conjecture}}
\newcommand{\EConjecture}{\end{Conjecture}}
\newcommand{\BProposition}{\begin{Proposition}}
\newcommand{\EProposition}{\end{Proposition}}
\newcommand{\BCorollary}{\begin{Corollary}}
\newcommand{\ECorollary}{\end{Corollary}}
\newcommand{\BKorollar}{\begin{Korollar}}
\newcommand{\EKorollar}{\end{Korollar}}
\newcommand{\BDefinition}{\begin{Definition}}
\newcommand{\EDefinition}{\end{Definition}}
\newcommand{\BLemma}{\begin{Lemma}}
\newcommand{\ELemma}{\end{Lemma}}
\newcommand{\BClaim}{\begin{Claim}}
\newcommand{\EClaim}{\end{Claim}}

\newcommand{\EExt}{{\operatorname{Ext}}}

\newcommand{\HHom}{{\operatorname{Hom}}}

\newcommand{\ra}{\rightarrow}

\newcommand{\la}{\lambda}

\numberwithin{equation}{subsection}
\begin{document}
\footnote{Supported
by the TMR-Network Algebraic Lie Theory ERB FMRX-CT97-0100 and
EPSRC Grant M22536 }

\title{Appendix}
\author{Steen Ryom-Hansen}
\address{Department of Mathematics\\City University\\
Northampton Square, London, EC1V 0HB, United Kingdom }
\maketitle 


\section{}

We shall follow the notation introduced in the paper of A.E. Parker [P].
It is then well known that for any finite dimensional $G$-module $M$ we 
have that 
\begin{equation} \label{5}
(M,\nabla (\mu))=\sum_{i} (-1)^{i} \dim \EExt^{i}_{G} (M,\nabla(\mu))
\end{equation}
Now the Lusztig conjecture (which is still unproved) 
can be formulated as a statement about 
$$ \mbox{ dim Ext}_G^i ( L(\la),\nabla(\mu)) $$
for $ \la $ and $ \mu $ in a certain region of the dominant weights.

\medskip
Especially the 
Lusztig conjecture would imply that for certain $ w, y \in W_p $ 
$$ \mbox{ dim Ext}_G^{l(w)-l(y)} ( L(w.\la),\nabla(y.\mu))=1 
$$

Our aim is to prove this statement for all $ w, y \in W_p $, i.e. even 
for weights in the region where the Lusztig conjecture is known not to hold.

\medskip
Our main methods to obtain this result will be the exact sequences from 
Corollary 3.4 together with a result of Jantzen on the translation 
functors applied to simple modules. It is a pleasure to thank the referee for many useful 
suggestions and especially for pointing out this result, which replaces our original 
argument. I also wish to thank H. H. Andersen for useful discussion related to
this work.

\section{Translation arguments}
Recall the translation functors $ T^{\mu}_{\la} $ and 
$ T^{\la}_{\mu} $  from section 3 of the paper.
Choose $ \la $ a regular, and $ \mu $ a semiregular weight. Choose moreover
$ w \in W_p $ such that $ ws < w $ for a simple reflection $s$.
It is then an easy, well known, fact that $ T^{\la}_{\mu} L(w.\mu) $ 
has simple head and socle, both isomorphic to $ L(ws.\la) $. Following 
[A] one can thus define modules $ U(w.\la) $ and $ Q(w.\la) $ by the exact 
sequences
\begin{equation}
0 \ra L(ws.\la) \ra T^{\la}_{\mu} L(w.\mu) \ra Q(w.\la) \ra 0
\end{equation}
\begin{equation}
0\ra U(w.\la) \ra Q(w.\la) \ra L(ws.\la) \ra 0
\end{equation}
It should be noticed that the action of the affine Weyl group on the weights 
is here the one of the above paper, rather than the one of [A].

\medskip

The following result is well known, see e.g. [J1], II.6.20.
\begin{Lemma}
For $w,y \in W_{p}$ such that $w.\la \in X(T)^{+}$ we have for $i
> l(w) - l(y)$ that
$$\EExt^{i}_{G} (L(w.\la),\nabla(y.\la))=0$$
\end{Lemma}

\subsection{}
The next two Lemmas relate the $\EExt$ groups we are
interested in to an $\EExt$ group involving $U(w.\lambda)$.
\begin{Lemma}
Assume $ys <y$. Then the following holds
$$\begin{array}{rcl}
\EExt^{l(w)-l(y)}_{G}(U(w.\lambda),\nabla(y.\lambda))&\cong&
\EExt_{G}^{l(w)-l(y)} (L(ws.\lambda),\nabla(ys.\lambda)) \text{ if
} ys.\lambda \in X(T)^{+}\\
\EExt_{G}^{l(w)-l(y)}(U(w.\lambda),\nabla(y.\lambda)) &\cong &
\EExt_{G}^{l(w)-l(y)-1} (L(ws.\lambda), \nabla(y.\lambda))\text{ if
} ys.\lambda \not\in X(T)^{+}
\end{array}$$
\end{Lemma}
\begin{pf}
By Lemma 2.6 of [A] we have that
$$\begin{array}{rcl}
\EExt^{i}_{G}(Q(w.\lambda), \nabla(y.\lambda)) &\cong &
\EExt^{i}_{G}(L(ws.\lambda)),\nabla(ys.\lambda)) \text{ if }
ys.\lambda \in X(T)^{+}\\
\EExt^{i}_{G} (Q(w.\lambda),\nabla(y.\lambda))&\cong &
\EExt^{i-1}_{G} (L(ws.\lambda),\nabla(y.\lambda)) \text{ if }
ys.\lambda \not\in X(T)^{+}
\end{array}$$
We insert this information into the long exact sequence that arises
from the application of $\HHom_{G}(-,\nabla(y.\lambda))$ to
(2.1.7). If $ys.\lambda \in X(T)^{+}$, part of the resulting
sequence is
$$\ra \EExt_{G}^{l(w)-l(y)} (L(ws.\lambda),\nabla(y.\lambda)) \ra
\EExt_{G}^{l(w)-l(y)} (L(ws.\lambda),\nabla(ys.\lambda))\ra$$
$$
\EExt_{G}^{l(w)-l(y)}(U(w.\lambda),\nabla(y.\lambda)) \ra
\EExt_{G}^{l(w)-l(y)+1}(L(ws.\lambda),\nabla(y.\lambda)) \ra$$
while for $ys.\lambda \not\in X(T)^{+}$ part of the resulting
sequence is
$$\ra \EExt_{G}^{l(w)-l(y)}(L(ws.\lambda),\nabla(y.\lambda))\ra
\EExt_{G}^{l(w)-l(w)-1} (L(ws.\lambda),\nabla(y.\lambda))\ra $$
$$\EExt_{G}^{l(w)-l(y)}(U(w.\lambda),\nabla(y.\lambda))\ra
\EExt_{G}^{l(w)-l(y)+1} (L(ws.\lambda),\nabla(y.\lambda))\ra$$
But by Lemma 2.1 the first and last terms in the two sequences are
zero; the Lemma is proved.
\end{pf}
\begin{Lemma}
Assume $sy > y$. Then the following holds
$$\EExt^{l(w)-l(y)}_{G}(U(w.\lambda), \nabla(y.\lambda))\cong
\EExt_{G}^{l(w)-l(y)-1} (L(ws.\lambda),\nabla(y.\lambda))$$
\end{Lemma}
\begin{pf}
This follows from the 
last Lemma and the exact sequence of Corollary 3 of the
paper, together with the fact that $ T^{\mu}_{\la} \, U(w. \la ) = 0 $.
\end{pf}
\subsection{}
After these preparatory Lemmas we can prove our main result.
\begin{Theorem}
For all $y,w \in W_{p}$ such that $y\leq w$ and $y.\lambda$,
$w.\lambda \in X(T)^{+}$ we have
$$\dim \EExt_{G}^{l(w)-l(y)}(L(w.\lambda),\nabla(y.\lambda))=1$$
\end{Theorem}
\begin{pf}
We proceed by induction on $l(w)$. T $ w=1 \Rightarrow y=1 $ and 
$$\dim \HHom_{G}(L(\lambda), \nabla(\lambda))=1$$
We then assume the result for $w^{\prime}$ with $l(w^{\prime}) <
l(w)$ and choose $s$ with $ws < w$ and $ws.\lambda \in X(T)^{+}$.
Then the Theorem holds for $ws$ and we get from Lemma 2.2 and
Lemma 2.3 that
$$\dim\EExt_{G}^{l(w)-l(y)} (U(w.\lambda),\nabla(y.\lambda))=1$$
So the Theorem would be a consequence of the isomorphism
\begin{equation}
\EExt_{G}^{l(w)-l(y)}(U(w.\lambda),\nabla(y.\lambda))\cong
\EExt_{G}^{l(w)-l(y)}(L(w.\lambda),\nabla(y.\lambda))
\end{equation}
We now claim that
\begin{equation}
[U(w.\lambda),L(z.\lambda)] \neq 0, \,z \neq w \Rightarrow l(w) -
l(z) \geq 2
\end{equation}
Believing this we would have from Lemma 2.1 that
\begin{equation}
\EExt_{G}^{l(w)-l(y)-k} (L(z.\lambda),\nabla(y.\lambda) )=0 \text{ for } k=0,1
\end{equation}
Now it is easy to see, (proof of prop. 2.8 (ii) of [A]) that
\begin{equation}
[U(w.\lambda),L(w.\lambda)]=1
\end{equation}
And then (2.4.1) would follow by filtering $U(w.\lambda)$ and
considering the terms of index $l(w)-l(y)$ and $l(w)-l(y) -1$ in
the long exact sequence given by the application of
$\HHom_{G}(-,\nabla(y.\lambda))$. So we must prove (2.2.2).

\medskip

But this follows from a result of Jantzen, part (b) of the proposition on 
page 299 in [J2], saying that 
\begin{equation}
[U(w.\lambda),L(z.\lambda)] 
\leq 2 \,[ \nabla( ws.\la ), L(z.\lambda)],\,\,\,\,\, \,z \neq w 
\end{equation}

\end{pf}

\end{document}